\documentclass{article}
\usepackage{amsmath,times}
\usepackage{amsthm}
\usepackage{amssymb}
\usepackage{latexsym}
\usepackage{amscd}
\usepackage[english]{babel}
\usepackage[latin1]{inputenc}
\newtheorem{theorem}{Theorem}[section]
\newtheorem{lemma}[theorem]{Lemma}
\newtheorem{proposition}[theorem]{Proposition}
\newtheorem{definition}[theorem]{Definition}
\newtheorem{corollary}[theorem]{Corollary}

\newtheorem{remark}[theorem]{Remark}
\begin{document} 

\title{Test configurations and Okounkov bodies}
\author{DAVID WITT NYSTRÖM}
\maketitle

\begin{abstract}
We associate to a test configuration of an ample line bundle a filtration of the section ring of the line bundle. Using the recent work of Boucksom-Chen we get a concave function on the Okounkov body whose law with respect to Lebesgue measure determines the asymptotic distribution of the weights of the test configuration. We show that this is a generalization of a well-known result in toric geometry. As an application, we prove that the pushforward of the Lebesgue measure on the Okounkov body is equal to a Duistermaat-Heckman measure of a certain deformation of the manifold. Via the Duisteraat-Heckman formula, we get as a corollary that in the special case of an effective $\mathbb{C^{\times}}$-action on the manifold lifting to the line bundle, the pushforward of the Lebesgue measure on the Okounkov body is piecewise polynomial.  
\end{abstract}

\tableofcontents

\section{Introduction}

\subsection{Okounkov bodies}

In \cite{Okounkov} Okounkov introduced a way to associate a convex body in $\mathbb{R}^n$ to any ample divisor on a $n$-dimensional projective variety. This procedure was later shown to work in a more general setting by Lazarsfeld-Musta\c{t}\u{a} in \cite{Lazarsfeld} and by Kaveh-Khovanskii in \cite{Kaveh} and \cite{Kaveh2}. 

Let $L$ be a big line bundle on a complex projective manifold $X$ of dimension $n.$ The Okounkov body of $L,$ denoted by $\Delta(L),$ is a convex subset of $\mathbb{R}^n,$ constructed in such a way so that the set-valued mapping $$\Delta: L \longmapsto \Delta(L)$$ has some very nice properties (for the explicit construction see Section \ref{sec1}). It is homogeneous, i.e. for any $k\in \mathbb{N}$ $$\Delta(kL)=k\Delta(L).$$ Here $kL$ denotes the the $k$:th tensor power of the line bundle $L.$ Secondly, the mapping is convex, in the sense that for any big line bundles $L$ and $L',$ and any $k,m\in \mathbb{N},$ the following holds $$\Delta(kL+mL')\supseteq k\Delta(L)+m\Delta(L'),$$ where the plus sign on the right hand side refers to Minkowski addition, i.e. $$A+B:=\{x+y : x\in A, y\in B\}.$$ Recall that the volume of a line bundle $L,$ denoted by $\textrm{vol}(L),$ is defined by $$\textrm{vol}(L):=\limsup_{k \to \infty}\frac{\dim H^0(kL)}{k^n/n!}.$$ By definition $L$ is big if $\textrm{vol}(L)>0.$ The third and crucial property, which makes Okounkov bodies useful as a tool in birational geometry, is that for any $L$ $$\textrm{vol}(L)=n!\textrm{vol}_{\mathbb{R}^n}(\Delta(L)).$$ where the volume of the Okounkov body is measured with respect to the standard Lesbesgue measure on $\mathbb{R}^n.$

\subsection{Test configurations}

Given an ample line bundle $L$ on $X,$ a class of algebraic deformations of the pair $(X,L),$ called test configurations, were introduced by Donaldson in \cite{Donaldson2}, generalizing a previous notion of Tian \cite{Tian} in the context of Fano manifolds. In short, a test configuration consists of: 
\begin{enumerate}
\item[(i)] a scheme $\mathcal{X}$ with a $\mathbb{C}^{\times}$-action $\rho,$
\item[(ii)] an $\mathbb{C}^{\times}$-equivariant line bundle $\mathcal{L}$ over $\mathcal{X},$ 
\item[(iii)] and a flat $\mathbb{C}^{\times}$-equivariant projection $\pi: \mathcal{X} \to \mathbb{C}$ such that $\mathcal{L}$ restricted to the fiber over $1$ is isomorphic to $L.$
\end{enumerate}
To a test configuration $\mathcal{T}$ there are associated discrete weight measures $\tilde{\mu}(\mathcal{T},k)$ (see Section \ref{sec3} for the definition). The asymptotics of the first moments of these measures are used to formulate stability conditions, such as $K$-stability, on the pair $(X,L).$ These  conditions are conjectured to be equivalent to the existence of a constant scalar curvature metric with Kähler form in $c_1(L),$ a conjecture which is sometimes called the Yau-Tian-Donaldson conjecture. This is one of the big open problems in Kähler geometry. By the works of e.g. Yau, Tian and Donaldson, a lot of progress has been made, in particular in the case of Kähler-Einstein metrics, i.e. when $L$ is a multiple of the canonical bundle. For more on this, we refer the reader to the expository article \cite{Sturm} by Phong-Sturm.  

When $L$ is assumed to be a toric line bundle on a toric variety with associated polytope $P,$ it was shown by Donaldson in \cite{Donaldson} that a test configuration is equivalent to specifying a concave rationally piecewise affine function on the polytope $P.$ This has made it possible to translate algebraic stability conditions on $L$ into geometric conditions on $P,$ which has proved very useful.

Heuristically, the relationship between a general line bundle $L$ and its Okounkov body is supposed to mimic the relationship between a toric line bundle and its associated polytope. Therefore, one would hope that one could translate a general test configuration into some geometric data on the Okounkov body. The main goal of this article is to show that this in fact can be done, thus presenting a generalization of the well-known toric picture referred to above, and described in greater detail in Section \ref{sec5}.

\subsection{The concave transform of a test configuration}

By a filtration $\mathcal{F}$ of the section ring $\oplus_k H^0(kL)$ we mean a vector space-valued map from $\mathbb{R}\times \mathbb{N},$ $$\mathcal{F}: (t,k) \longmapsto \mathcal{F}_t H^0(kL),$$ such that for any $k,$ $\mathcal{F}_t H^0(kL)$ is a family of subspaces of $H^0(kL)$ that is decreasing and left-continuous in $t.$ $\mathcal{F}$ is said to be multiplicative if $$ (\mathcal{F}_t H^0(kL))( \mathcal{F}_s H^0(mL))\subseteq  \mathcal{F}_{t+s} H^0((k+m)L),$$ it is left-bounded if for all $k$ $$\mathcal{F}_{-t} H^0(kL)=H^0(kL) \qquad{} \textrm{for} \qquad{} t\gg 1,$$ and is said to linearly right-bounded if there exist a $C$ such that $$\mathcal{F}_t H^0(kL)=\{0\}\qquad{} \textrm{for} \qquad{}t\geq Ck.$$ The filtration $\mathcal{F}$ is called admissible if it has all the above properties.

Given a filtration $\mathcal{F}$, one may associate discrete measures $\nu(\mathcal{F},k)$ on $\mathbb{R}$ in the following way $$\nu(\mathcal{F},k):=\frac{1}{k^n}\frac{d}{dt}(-\dim \mathcal{F}_{tk} H^0(kL)),$$ where the differentiation is done in the sense of distributions.

In a recent preprint \cite{Chen} Boucksom-Chen show how any admissible filtration $\mathcal{F}$ of the section ring $\oplus_k H^0(kL)$ of $L$ gives rise to a concave function $G[\mathcal{F}]$ on the Okounkov body $\Delta(L)$ of $L.$ $G[\mathcal{F}]$ is called the concave transform of $\mathcal{F}$. The main result of  \cite{Chen}, Theorem A, states that the discrete measures $\nu(\mathcal{F},k)$ converge weakly as $k$ tends to infinity to $G[\mathcal{F}]_*d\lambda_{|\Delta(L)},$ the push-forward of the Lebesgue measure on $\Delta(L)$ with respect to the concave transform of $\mathcal{F}$.

Let $\mathcal{T}$ be a test configuration on $(X,L).$ Given a section $s\in H^0(kL),$ there is a unique invariant meromorphic extension to configuration scheme $\mathcal{X}.$ Using the vanishing order of this extension along the central fiber of $\mathcal{X}$ we define a filtration of the section ring $\oplus_k H^0(kL)$, which we show has the property that for any $k$ $$\tilde{\mu}(\mathcal{T},k)=\nu(\mathcal{F},k).$$  We will denote the associated concave transform by $G[\mathcal{T}]$. Combined with Theorem A of \cite{Chen} we thus get our first main result.

\begin{theorem} \label{thmthm1}
Given a test configuration $\mathcal{T}$ of $L$ there is a concave function $G[\mathcal{T}]$ on the Okounkov body $\Delta(L)$ such that the measures $\tilde{\mu}(\mathcal{T},k)$ converge weakly as $k$ tends to infinity to the measure $G[\mathcal{T}]_*d\lambda_{|\Delta(L)}.$
\end{theorem} 

We embed our test configuration into $\mathbb{C}$ times a projective space $\mathbb{P}^N,$ so that the associated action comes from a $\mathbb{C}^{\times}$-action on $\mathbb{P}^N.$ This we can always do (see e.g. \cite{Ross2}). The manifold $X$ lies embedded in $\mathbb{P}^N,$ and we thus via the action get a family $X_{\tau}$ of subamnifolds. As $\tau$ tends to $0,$ $X_{\tau}$ converges in the sense of currents to an algebraic cycle $|X_0|$ (see \cite{Donaldson}). We let $\omega_{FS}$ denote the Fubini-Study on $\mathbb{P}^N.$ Restricted to $X_{\tau}$ the $(n,n)$-form $\omega_{FS}^n/n!$ defines a positive measure, that as $\tau$ goes to zero converges to a positive measure $d\mu_{FS},$ the Fubini-Study volume form on $|X_0|.$  There is also a Hamiltonian function $h$ for the $S^1$-action. Using a result of Donaldson in \cite{Donaldson} and Theorem \ref{thmthm1} we can relate this picture with the concave transform by the following Corollary.      

\begin{corollary} \label{retta}
Assume that we have embedded the test configuration $\mathcal{T}$ in some $\mathbb{P}^N\times \mathbb{C},$ let $h$ denote the corresponding Hamiltonian and $d\mu_{FS}$ the positive measure on $|X_0|$ defined above. Then we have that $$h_* d\mu_{FS}=G[\mathcal{T}]_*d\lambda_{|\Delta(L)}.$$ 
\end{corollary}

If $|X_0|$ is a smooth manifold, on which the $S^1$-action is effective, the measure $h_* d\mu_{FS}$ is the sort of measure studied by Duistermaat-Heckman in \cite{Duistermaat}. They prove that such a Duistermaat-Heckman measure is piecewise polynomial, i.e. the distribution function with respect to Lebesgue measure on $\mathbb{R}$ is piecewise polynomial. For a product test configuration, $|X_0|\cong X,$ therefore we can apply the result of Duistermaat-Heckman to get the following.  

\begin{corollary} \label{riwte}
Assume that there is a $\mathcal{C}^{\times}$-action on $X$ which lifts to $L,$ and that the corresponding $S^1$-action is effective. If we denote the associated product test configuration by $\mathcal{T},$ the concave transform $G[\mathcal{T}]$ is such that the pushforward measure $G[\mathcal{T}]_*d\lambda_{|\Delta(L)}$ is piecewise polynomial.
\end{corollary} 

We also consider the case of a product test configuration, which means that there is an algebraic $\mathbb{C}^{\times}$-action $\rho$ on the pair $(X,L).$ We let $\varphi$ be a positive $S^1$-invariant metric on $L.$ Using the action $\rho,$ we get a geodesic ray $\varphi_t$ of positive metrics on $L$ such that $\varphi_1=\varphi.$ Let us denote the $t$ derivative at the point one by $\dot{\varphi}.$ It is a real-valued function on $X.$ There is also a natural volume element, given by $dV_{\varphi}:=(dd^c\varphi)^n/n!.$ By the function $\dot{\varphi}/2$ we can push forward the measure $dV_{\varphi}$ to a measure on $\mathbb{R},$ which we denote by $\mu_{\varphi}.$ This measure does not depend on the particular choice of positive $S^1$-invariant metric $\varphi.$ In fact, we have the following.

\begin{theorem} \label{thmthm2}
If we denote the product test configuration by $\mathcal{T},$ and the corresponding concave transform by $G[\mathcal{T}]$, then for any positive $S^1$-invariant metric $\varphi$ it holds that $$\mu_{\varphi}=G[\mathcal{T}]_*d\lambda_{|\Delta(L)}.$$
\end{theorem}  

The proof uses Theorem \ref{thmthm1} combined with the approach of Berndtsson in \cite{Berndtsson}, but is simpler in nature.

Phong-Sturm have in their articles \cite{Sturm} and \cite{Sturm2} shown that the pair of a test configuration $\mathcal{T}$ and a positive metric $\varphi$ on $L$ canonically determines a $C^{1,1}$ geodesic ray of positive metrics on $L$ emanating from $\varphi.$ We conjecture that the analogue of Theorem \ref{thmthm2} is true also in that more general case.

%From the symplectic viewpoint, $\dot{\varphi}/2$ is nothing but the moment map of $X$ related to the action $\rho.$ Thus using the Duistermaat-Heckman formula, see \cite{Duistermaat}, we get as a corollary the following.

%\begin{corollary} \label{crlcrl2}
%Assume that there is an algebraic $\mathbb{C}^{\times}$-action $\rho$ on the pair $(X,L).$ Then there exists a concave function $G[\mathcal{T}]$ on the Okounkov body $\Delta(L),$ namely the concave transform of the corresponding test configuration, such that the pushforward measure $$G[\mathcal{T}]_*d\lambda_{|\Delta(L)}$$ is piecewise polynomial of degree less than $n$.  
%\end{corollary}   

\subsection{Organization of the paper}

The definition of Okounkov bodies and some fundamental results concerning them is in Section \ref{sec1}, using \cite{Lazarsfeld} by Lazarsfeld-Musta\c{t}\u{a} as our main reference.

Section \ref{sec2} is devoted to describing the setup, definitions and main results of the article \cite{Chen} by Boucksom-Chen on the concave transform of filtrations.

Section \ref{sec3} contains a brief introduction to test configurations, following mainly Donaldson in \cite{Donaldson2} and \cite{Donaldson}.

We discuss embeddings of test configurations in Section \ref{secny}, and link it to certain Duistermaat-Heckman measures.

In Section \ref{sec4} we show how to construct the associated filtration to a test configuration, and prove Theorem \ref{thmthm1}, Corollary \ref{retta} and Corollary \ref{riwte}.

Section \ref{sec5} concerns toric test configurations. We show that what we have done is a generalization of the toric picture, by proving that in the toric case, the concave transform is identical to the function on the polytope considered by Donaldson in \cite{Donaldson2}.

Relying on the work of Ross-Thomas in \cite{Ross} and \cite{Ross2}, we obtain in Section \ref{sec6} an explicit description of the concave transforms corresponding to a special class of test configurations, namely those arising from a deformation to the normal cone with respect to some subscheme.

In Section \ref{sec7} we study the case of product test configurations, and relate it to geodesic rays of positive hermitian metrics. Hence we prove Theorem \ref{thmthm2}.

\subsection{Acknowledgements}

We wish to thank Robert Berman, Bo Berndtsson, Sébastien Boucksom, Julius Ross and Xiaowei Wang for many interesting discussions relating to the topic of this paper.

\section{The Okounkov body of a line bundle} \label{sec1}

Let $\Gamma$ be a subset of $\mathbb{N}^{n+1}$, and suppose that it is a semigroup with respect to vector addition, i.e. if $\alpha$ and  $\beta$ lie in $\Gamma,$ then the sum $\alpha+\beta$ should also lie in $\Gamma$. We denote by $\Sigma(\Gamma)$ the closed convex cone in $\mathbb{R}^{n+1}$ spanned by $\Gamma.$ 

\begin{definition}
The Okounkov body $\Delta(\Gamma)$ of $\Gamma$ is defined by $$\Delta(\Gamma):=\{\alpha: (\alpha,1)\in \Sigma(\Gamma)\}\subseteq  \mathbb{R}^{n}.$$
\end{definition}

Since by definition $\Sigma(\Gamma)$ is convex, and any slice of a convex body is itself convex, it follows that the Okounkov body $\Delta(\Gamma)$ is convex.

By $\Delta_k(\Gamma)$ we will denote the set $$\Delta_k(\Gamma):=\{\alpha: (k\alpha,k)\in \Gamma\}\subseteq \mathbb{R}^{n}.$$ It is clear that for all non-negative $k,$ $$\Delta_k(\Gamma)\subseteq \Delta(\Gamma)\cap((1/k)\mathbb{Z})^n.$$

We will explain the procedure, which is due to Okounkov (see \cite{Okounkov}), of associating a semigroup to a big line bundle.

Let $X$ be a complex compact projective manifold of dimension $n$, and $L$ a holomorphic line bundle, which we will assume to be big. Suppose we have chosen a point $p$ in $X,$ and local holomorphic coordinates $z_1,...,z_n$ centered at $p$, and let $e_p\in H^0(U,L)$ be a local trivialization  of $L$ around $p.$ If we divide a section $s\in H^0(X,kL)$ by $e_p^k$ we get a local holomorphic function. It has an unique represention as a convergent power series in the variables $z_i,$ $$\frac{s}{e_p^k}=\sum a_{\alpha}z^{\alpha},$$ which for convenience we will simply write as $$s=\sum a_{\alpha}z^{\alpha}.$$ We consider the lexicographic order on the multiindices $\alpha$, and let $v(s)$ denote the smallest index $\alpha$ such that $a_{\alpha} \neq 0.$

\begin{definition}
Let $\Gamma(L)$ denote the set $$\{(v(s),k): s\in H^0(kL), k\in \mathbb{N}\} \subseteq \mathbb{N}^{n+1}.$$ It is a semigroup, since for $s\in H^0(kL)$ and $t\in H^0(mL)$ 
\begin{equation*} 
v(st)=v(s)+v(t).
\end{equation*}
The Okounkov body of $L$, denoted by $\Delta(L)$, is defined as the Okounkov body of the associated semigroup $\Gamma(L).$
\end{definition}

We write $\Delta_k(\Gamma(L))$ simply as $\Delta_k(L).$

\begin{remark}
Note that the Okounkov body $\Delta(L)$ of a line bundle $L$ in fact depends on the choice of point $p$ in $X$ and local coordinates $z_i.$ We will however supress this in the notation, writing $\Delta(L)$ instead of the perhaps more proper but cumbersome $\Delta(L,p,(z_i)).$
\end{remark}

From the article \cite{Lazarsfeld} by Lazarsfeld-Musta\c{t}\u{a} we recall some results on Okounkov bodies of line bundles.

\begin{lemma} \label{points}
The number of points in $\Delta_k(L)$ is equal to the dimension of the vector space $H^0(kL).$
\end{lemma}

\begin{lemma}
We have that $$\Delta(L)=\overline{\cup_{k=1}^{\infty}\Delta_k(L)}.$$
\end{lemma}

\begin{lemma} \label{okounkbound}
The Okounkov body $\Delta(L)$ of a big line bundle is a bounded hence compact convex body.
\end{lemma}

\begin{definition}
The volume of a line bundle $L,$ denoted by $\textrm{vol}(L),$ is defined by $$\textrm{vol}(L):=\limsup_{k \to \infty}\frac{\dim H^0(kL)}{k^n/n!}.$$ 
\end{definition}

The most important property of the Okounkov body is its relation to the volume of the line bundle, described in the following theorem.

\begin{theorem} \label{Lathm}
For any big line bundle it holds that $$\textrm{vol}(L)=n!\textrm{vol}_{\mathbb{R}^n}(\Delta(L)),$$ where the volume of the Okounkov body is measured with respect to the standard Lesbesgue measure on $\mathbb{R}^n.$
\end{theorem}

For the proof see \cite{Lazarsfeld}.

\section{The concave transform of a filtered linear series} \label{sec2}

In this section, we will follow Boucksom-Chen in \cite{Chen}.

First we recall what is meant by a filtration of a graded algebra.

\begin{definition}
By a filtration $\mathcal{F}$ of a graded algebra $\oplus_k V_k$ we mean a vector space-valued map from $\mathbb{R}\times \mathbb{N},$ $$\mathcal{F}: (t,k) \longmapsto \mathcal{F}_t V_k,$$ such that for any $k,$ $\mathcal{F}_t V_k$ is a family of subspaces of $V_k$ that is decreasing and left-continuous in $t.$ 
\end{definition}

In \cite{Chen} Boucksom-Chen consider certain filtrations which behaves well with respect to the multiplicative structure of the algebra.

They give the following definition. 

\begin{definition}
Let $\mathcal{F}$ be a filtration of a graded algebra $\oplus_k V_k$. We shall say that 
\begin{enumerate}
\item[(i)] $\mathcal{F}$ is multiplicative if  $$(\mathcal{F}_t V_k)(\mathcal{F}_s V_m)\subseteq \mathcal{F}_{t+s}V_{k+m}$$
for all $k,m \in \mathbb{N}$ and $s,t \in \mathbb{R}$.
\item[(ii)] $\mathcal{F}$ is pointwise left-bounded if for each $k$ $\mathcal{F}_t V_k=V_k$ for some $t$.
\item[(iii)] $\mathcal{F}$ is linearly right-bounded if there exist a constant $C$ such that for all $k,$ $\mathcal{F}_{kC}V_k=\{0\}.$
\end{enumerate}

A filtration $\mathcal{F}$ is said to be admissible if it is multiplicative, pointwise left-bounded and linearly right-bounded.
\end{definition}

Given a line bundle $L$ on $X,$ its section ring $\oplus_k H^0(kL)$ is a graded algebra.

Boucksom-Chen in \cite{Chen} show how an admissible filtration on the section ring $\oplus_k H^0(kL)$ of a big line bundle $L$ gives rise to a concave function on the Okounkov body $\Delta(L)$. We will review how this is done.

First let us define the following set $$\Delta_{k,t}(L,\mathcal{F}):=\{v(s)/k: s\in \mathcal{F}_t H^0(kL)\}\subseteq \mathbb{R}^n,$$ where as before $v(s)=\alpha$ if locally $$s=Cz^{\alpha}+ \textrm{higher order terms},$$ $C$ being some nonzero constant. From the definition it is clear that $$\Delta_{k,t}(L,\mathcal{F})\subseteq \Delta_k(L),$$ since $$\Delta_k(L)=\{v(s)/k: s\in H^0(kL)\}$$ and $\mathcal{F}_t H^0(kL)\subseteq H^0(kL).$ Similarly as in Lemma \ref{points}, from \cite{Lazarsfeld} we get that 
\begin{equation} \label{dimpoint}
|\Delta_{k,t}(L,\mathcal{F})|=\dim \mathcal{F}_t H^0(kL),
\end{equation}
where $|.|$ denotes the cardinality of the set.

For each $k$ we may define a function $G_k$ on $\Delta_k(L)$ by letting $$G_k(\alpha):=\sup \{t: \alpha\in \Delta_{k,t}(L,\mathcal{F})\}.$$ From the assumption that $\mathcal{F}$ is both left- and right-bounded it follows that $G_k$ is well-defined and real-valued. 

\begin{lemma} \label{rara}
If we denote by $\nu_k(L)$ the sum of dirac measures at the points of $\Delta_k(L),$ i.e. $$\nu_k(L):=\sum_{\alpha\in \Delta_k(L)}\delta_{\alpha},$$ then we have that $${G_k}_*\nu_k(L)=\frac{d}{dt}(-\dim \mathcal{F}_t H^0(kL)).$$
\end{lemma}

\begin{proof}
By equation (\ref{dimpoint}) and the definition of $G_k$ we have that 
\begin{equation} \label{dimeq1}
\dim \mathcal{F}_t H^0(kL)=|\Delta_{k,t}(L,\mathcal{F})|=\int_{\{G_k\geq t\}}d\nu_k(L)=\int_t^{\infty}(G_k)_*(\nu_k(L)).
\end{equation}
The lemma now follows by differentiating the equation (\ref{dimeq1}).
\end{proof}

On the union $\cup_{k=1}^{\infty}\Delta_k(L)$ one may define the function $$G[\mathcal{F}](\alpha):=\sup\{G_k(\alpha)/k: \alpha\in \Delta_k(L)\}.$$ By Boucksom-Chen in \cite{Chen}, or Witt Nyström in \cite{Witt}, one then gets that the function $G[\mathcal{F}]$ extends to a concave and therefore continuous function on the interior of $\Delta(L).$ In fact one gets that $G[\mathcal{F}]$ is not only the supremum but also the limit of $G_k/k,$ i.e. for any $p\in \Delta(L)^{\circ}$ $$G[\mathcal{F}](p)=\lim_{k\to \infty}G_k(\alpha_k)/k,$$ for any sequence $\alpha_k$ converging to $p.$  

\begin{remark}
To show how this fits into the framework of \cite{Witt}, we note that if we let $$\tilde{G}(\alpha,k):=G_k(\alpha/k),$$ then $\tilde{G}$ is a function on $\Gamma(L).$ By the multiplicity of $\mathcal{F}$ it follows that $\tilde{G}$ is superadditive, and by the linear right-boundedness, $\tilde{G}$ is going to be linearly bounded from above. Thus one may apply the results of \cite{Witt} to this function. 
\end{remark}

The main result of \cite{Chen}, Theorem A, is that we also have weak convergence of measures.

\begin{theorem} \label{sebthm}
The measures $$\frac{1}{k^n}({(G_k/k)}_*\nu_k(L))$$ converge weakly to the measure $$G[\mathcal{F}]_*d\lambda_{|\Delta(L)}$$ as $k$ tends to infinity, where $d\lambda_{|\Delta(L)}$ denotes the Lebesgue measure on $\mathbb{R}^n$ restricted to $\Delta(L).$
\end{theorem}

\section{Test configurations} \label{sec3}

We will give a very brief introduction to the subject of test configurations. Our main references are the articles \cite{Donaldson2} and \cite{Donaldson} by Donaldson.

First the definition of a test configuration, as introduced by Donaldson in \cite{Donaldson2}.

\begin{definition}
A test configuration $\mathcal{T}$ for an ample line bundle $L$ over $X$ consists of:
\begin{enumerate}
\item[(i)] a scheme $\mathcal{X}$ with a $\mathbb{C}^{\times}$-action $\rho,$
\item[(ii)] an $\mathbb{C}^{\times}$-equivariant line bundle $\mathcal{L}$ over $\mathcal{X},$ 
\item[(iii)] and a flat $\mathbb{C}^{\times}$-equivariant projection $\pi: \mathcal{X} \to \mathbb{C}$ where $\mathbb{C}^{\times}$ acts on $\mathbb{C}$ by multiplication, such that $\mathcal{L}$ is relatively ample, and such that if we denote by $X_1:=\pi^{-1}(1)$, then $\mathcal{L}_{|X_1} \to X_1$ is isomorphic to $rL \to X$ for some $r>0.$ 
\end{enumerate}
\end{definition}

By rescaling we can for our purposes without loss of generality assume that $r=1$ in the definition. 

A test configuration is called a product test configuration if there is a $\mathbb{C}^{\times}$-action $\rho'$ on $L \to X$ such that $\mathcal{L}=L\times \mathbb{C}$ with $\rho$ acting on $L$ by $\rho'$ and on $\mathbb{C}$ by multiplication. A test configuration is called trivial if it is a product test configuration with the action $\rho'$ being the trivial $\mathbb{C}^{\times}$-action.

Since the zero-fiber $X_0:=\pi^{-1}(0)$ is invariant under the action $\rho$, we get an induced action on the space $H^0(kL_0),$ also denoted by $\rho,$ where we have denoted the restriction of $\mathcal{L}$ to $X_0$ by $L_0$. Specifically, we let $\rho(\tau)$ act on a section $s\in H^0(kL_0)$ by 
\begin{equation} \label{graction}
(\rho(\tau)(s))(x):=\rho(\tau)(s(\rho^{-1}(\tau)(x))).
\end{equation}

\begin{remark}
Some authors refer to the inverted variant $$(\rho(\tau)(s))(x):=\rho^{-1}(\tau)(s(\rho(\tau)(x)))$$ as the induced action. This is only a matter of convention, but one has to be aware that all the weights as defined below changes sign when changing from one convention to the other.  
\end{remark}

Any vector space $V$ with a $\mathbb{C}^{\times}$-action can be split into weight spaces $V_{\eta_i}$ on which $\rho(\tau)$ acts as multiplication by $\tau^{\eta_i}$, (see e.g. \cite{Donaldson2}). The numbers $\eta_i$ with non-trivial weight spaces are called the weights of the action. Thus we may write $H^0(kL_0)$ as $$H^0(kL_0)=\oplus_{\eta}V_{\eta}$$ with respect to the induced action $\rho.$ 

In \cite{Sturm}, Lemma 4, Phong-Sturm give the following linear bound on the absolute value of the weights.

\begin{lemma} \label{phong}
Given a test configuration there is a constant $C$ such that $$|\eta_i|<Ck$$ whenever $\dim V_{\eta_i}>0$.
\end{lemma}

There is an associated weight measure on $\mathbb{R}:$ $$\mu(\mathcal{T},k):=\sum_{\eta=-\infty}^{\infty}\dim V_{\eta}\delta_{\eta},$$ and also the rescaled variant 
\begin{equation} \label{renorm}
\tilde{\mu}(\mathcal{T},k):=\frac{1}{k^n}\sum_{\eta=-\infty}^{\infty}\dim V_{\eta}\delta_{k^{-1}\eta}.
\end{equation}

The first moment of the measure $\mu(\mathcal{T},k),$ which we will denote by $w_k,$ thus equals the sum of the weights $\eta_i$ with multiplicity $\dim V_{\eta_i}$. It can also be seen as the weight of the induced action on the top exterior power of $H^0(kL_0).$ The total mass of $\mu(\mathcal{T},k)$ is $\dim H^0(kL_0),$ which we will denote by $d_k.$ By the flatness of $\pi$ it follows that for $k$ large it will be equal to $\dim H^0(kL)$ (see e.g. \cite{Ross}). One is interested in the asymptotics of the weights, and from the equivariant Riemann-Roch theorem one gets that there is an asymptotic expansion in powers of $k$ of the expression $w_k/kd_k$ (see e.g. \cite{Donaldson2}), 

$$\frac{w_k}{kd_k}=F_0-k^{-1}F_1+O(k^{-2}).$$ $F_1$ is called the Futaki invariant of $\mathcal{T}$, and will be denoted by $F(\mathcal{T})$.

\begin{definition}
A line bundle $L$ is called K-semistable if for all test configurations $\mathcal{T}$ of $L$ over $X,$ it holds that $F(\mathcal{T})\geq 0$. $L$ is called K-stable if it is K-semistable and furthermore $F(\mathcal{T})=0$ iff  $\mathcal{T}$ is a product test configuration.
\end{definition}

Donaldson has conjectured that $L$ being $K$-stable is equivalent to the existence of a positive constant scalar curvature hermitian metric with Kähler form in $c_1(L)$ (see \cite{Donaldson2}, \cite{Donaldson} and the expository article \cite{Sturm3}).

\section{Embeddings of test configurations} \label{secny}

One way to construct a test configuration of a pair $(X,L)$ is by using a Kodaira embedding of $(X,L)$ into $(\mathbb{P}^N, \mathcal{O}(1))$ for some $N$. If $\rho$ is a $\mathbb{C}^{\times}$-action on $\mathbb{P}^N,$ this gives rise to a product test configuration of $(\mathbb{P}^N, \mathcal{O}(1)).$ If we restrict to the image of $\rho$'s action on $(X,L),$ we end up with a test configuration of $(X,L).$ A basic fact (see e.g. \cite{Ross2}) is that all test configurations arise this way, so that one may embed $\mathcal{X}$ into $\mathbb{P}^N\times \mathbb{C}$ for some $N,$ the action $\rho$ coming from a $\mathbb{C}^{\times}$-action on $\mathbb{P}^N.$ 

Let $\mathcal{T}$ be a test configuration, and assume that we have chosen an embeddding as above. Let $z_i$ be homogeneous coordinates on $\mathbb{P}^N,$ and let us define the following functions $$h_{ij}:=\frac{z_i \bar{z}_j}{||z||^2}.$$ We assume that we have chosen our coordinates so that the metric $||z||^2$ is invariant under the corresponding $S^1$-action on $\mathbb{C}^{N+1}$. Then the infinitesimal generator of the action $\rho$ is given by a hermitian matrix $A.$ We define a real-valued function $h$ on $\mathbb{P}^N$ by $$h:=\sum A_{ij}h_{ij}.$$ It is a Hamiltonian for the $S^1$-action (see \cite{Donaldson}). Let $\omega_{FS}$ denote the Fubini-Study form on $\mathbb{P}^N.$ The zero-fiber $X_0$ of the test configuration can via the embedding be identified with subsheme of $\mathbb{P}^N,$ invariant under the action of $\rho.$ By $|X_0|$ we will denote the corresponding algebraic cycle, and we let $[X_0]$ denote its integration current. The wedge product of $[X_0]$ with the positive $(n,n)$-form $\omega_{FS}^n/{n!}$ gives a positive measures, $d\mu_{FS},$ with $|X_0|$ as its support. We have the following proposition.

\begin{proposition} \label{reta}
In the setting as above, the normalized weight measures $\tilde{\mu}(\mathcal{T},k)$ of the test configuration converges weakly as $k$ tends to infinity to the pushforward of the measure $d\mu_{FS}$ with respect to the Hamiltonian $h,$ $$\tilde{\mu}(\mathcal{T},k) \rightarrow h_*d\mu_{FS}.$$
\end{proposition}

\begin{proof}
This is essentially just a reformulation of a result by Donaldson in \cite{Donaldson}. Using the weight measures $\tilde{\mu}(\mathcal{T},k)$, Equation (20) in the proof of Proposition 3 in \cite{Donaldson} says that $$\int_{\mathbb{R}} x^r d\tilde{\mu}(\mathcal{T},k)=\int _{|X_0|}h^r d\mu_{FS}+o(1).$$ for any positive integer $r.$  In other words, for all such $r,$ the $r$-moments of the measures $\tilde{\mu}(\mathcal{T},k)$ converge to the $r$-moment of the pushforward measure $h_*d\mu_{FS}.$ But then it is classical that this implies weak convergence of measures. 
\end{proof}

The measure $h_* d\mu_{FS}$ is the sort of measure studied by Duistermaat-Heckman in \cite{Duistermaat}. They consider a smooth symplectic manifold $M$ with symplectic form $\sigma,$ and an effective Hamiltonian torus action on $M.$ This gives rise to a moment mapping $J$, which is a map from $M$ to the dual of the Lie algebra of the torus, which we can naturally identify with $\mathbb{R}^k,$ $k$ being the dimension of the torus (we refer the reader to \cite{Duistermaat} for the definitions). There is a natural volume measure on $M,$ given by $\sigma^n/n!,$ called the Liouville measure. The pushforward of the Liouville measure with the moment map $J$, $J_*(\sigma^n/n!),$ is called a Duistermaat-Heckman measure. They prove that it is absolute continuous with respect to Lebesgue measure on $\mathbb{R}^k,$ and provide an explicit formula, in the literature referred to as the Duistermaat-Heckman formula, for the density function $f.$ As a corollary they get the following.

\begin{theorem} \label{DuiHeck}
The density function $f$ of the measure $J_*(\sigma^n/n!)$ is a polynomial of degree less that the dimension of $M$ on each connected component of the set of regular values of the moment map $J.$
\end{theorem}

In our setting the Liouville measure is given by $d\mu_{FS},$ and the moment map $J$ is simply given by the Hamiltonian $h.$ Thus when all components of the algebraic cycle $|X_0|$ are smooth manifolds, and the action is effective, we can apply Theorem \ref{DuiHeck} to our measure $h_*d\mu_{FS}$ and conclude that it is a piecewise polynomial measure on $\mathbb{R}.$ In general of course some components of $|X_0|$ may have singularities. However, one case where we know that $X_0$ is a smooth manifolld is when we have a product test configuration, because then $X_0=X.$ Hence we get the following.

\begin{proposition} \label{prodtest}
For a product test configuration, with a corresponding effective $S^1$-action, it holds that the law of the asymptotic distribution of its weights is piecewise polynomial.  
\end{proposition} 

\begin{proof}
By Proposition \ref{reta} the law of the asymptotic distribution of weights is given by the measure $h_*d\mu_{FS}$ and by the remarks above we can use Theorem \ref{DuiHeck} to conclude that $h_*d\mu_{FS}$ is piecewise polynomial.
\end{proof}

\section{The concave transform of a test configuration} \label{sec4}

Given a test configuration $\mathcal{T}$ of $L$ we will show how to get an associated filtration $\mathcal{F}$ of the section ring $\oplus_k H^0(kL)$. 

First note that the $\mathbb{C}^{\times}$-action $\rho$ on $\mathcal{L}$ via the equation (\ref{graction}) gives rise to an induced action on $H^0(\mathcal{X}, k\mathcal{L})$ as well as $H^0(\mathcal{X}\setminus X_0,k\mathcal{L}),$ since $\mathcal{X}\setminus X_0$ is invariant.

Let $s\in H^0(kL)$ be a holomorphic section. Then using the $\mathbb{C}^{\times}$-action $\rho$ we get a canonical extension $\bar{s}\in H^0(\mathcal{X}\setminus X_0,k\mathcal{L})$ which is invariant under the action $\rho$, simply by letting 
\begin{equation} \label{toroto}
\bar{s}(\rho(\tau)x):=\rho(\tau)s(x)
\end{equation} 
for any $\tau\in \mathbb{C}^{\times}$ and $x\in X.$ 

We identify the coordinate $t$ with the projection function $\pi(x),$ and we also consider it as a section of the trivial bundle over $\mathcal{X}.$ Exactly as for $H^0(\mathcal{X}, k\mathcal{L}),$ $\rho$ gives rise to an induced action on sections of the trivial bundle, using the same formula (\ref{graction}).  We get that
\begin{equation} \label{useful}
(\rho(\tau)t)(x)=\rho(\tau)(t(\rho^{-1}(\tau)x)=\rho(\tau)(\tau^{-1}t(x))=\tau^{-1}t(x),
\end{equation}
where we used that $\rho$ acts on the trivial bundle by multiplication on the $t$-coordinate. Thus $$\rho(\tau)t=\tau^{-1}t,$$ which shows that the section $t$ has weight $-1.$ 

By this it follows that for any section $s\in H^0(kL)$ and any integer $\eta,$ we get a section $t^{-\eta}\bar{s}\in H^0(\mathcal{X}\setminus X_0,k\mathcal{L}),$ which has weight $\eta.$ 

\begin{lemma} \label{igenlemma}
For any section $s\in H^0(kL)$ and any integer $\eta$ the section $t^{-\eta}\bar{s}$ extends to a meromorphic section of $k\mathcal{L}$ over the whole of $\mathcal{X},$ which we also will denote by $t^{-\eta}\bar{s}.$
\end{lemma}

\begin{proof}
It is equivalent to saying that for any section $s$ there exists an integer $\eta$ such that $t^{\eta}\bar{s}$ extends to a holomorphic section $S\in H^0(\mathcal{X},k\mathcal{L}).$ By flatness, which was assumed in the definition of a test configuration, the direct image bundle $\pi_*\mathcal{L}$ is in fact a vector bundle over $\mathbb{C}$. Thus it is trivial, since any vector bundle over $\mathbb{C}$ is trivial. Therefore there exists a global section $S'\in H^0(\mathcal{X},k\mathcal{L})$ such that $s=S'_{|X}.$ On the other hand, as for $H^0(kL_0),$ $H^0(\mathcal{X},k\mathcal{L})$ may be decomposed as a direct sum of invariant subspaces $W_{\eta'}$ such that $\rho(\tau)$ restricted to $W_{\eta'}$ acts as multiplication by $\tau^{\eta'}$. Let us write 
\begin{equation} \label{decomp}
S'=\sum S'_{\eta'},
\end{equation}
where $S_{\eta'}\in W_{\eta'}.$ Restricting the equation (\ref{decomp}) to $X$ gives a decomposition of $s,$ $$s=\sum s_{\eta'},$$ where $s_{\eta'}:={S'_{\eta'}}_{|X}.$ From (\ref{toroto}) and the fact that $S'_{\eta'}$ lies in $W_{\eta'}$ we get that for $x\in X$ and $\tau\in \mathbb{C}^{\times}$ we have that 
\begin{eqnarray*}
\bar{s}_{\eta'}(\rho(\tau)(x))=\rho(\tau)(s_{\eta'}(x))=\rho(\tau)(S'_{\eta'}(x))=(\rho(\tau)S'_{\eta'})(\rho(\tau)(x)))=\\ =\tau^{\eta'}S'_{\eta'}(\rho(\tau)(x)),
\end{eqnarray*} 
and therefore $\bar{s}_{\eta'}=\tau^{\eta'}S'_{\eta'}.$ Since trivially $$\bar{s}=\sum \bar{s}_{\eta'}$$ it follows that $t^{\eta}\bar{s}$ extends holomorphically as long as $\eta\geq \max -\eta'.$  
\end{proof}

\begin{definition}
Given a test configuration $\mathcal{T}$ we define a vector space-valued map $\mathcal{F}$ from $\mathbb{Z}\times \mathbb{N}$ by letting $$(\eta,k) \longmapsto \{s\in H^0(kL): t^{-\eta}\bar{s}\in H^0(\mathcal{X},k\mathcal{L})\}=:\mathcal{F}_{\eta}H^0(kL).$$
\end{definition}
It is immediate that $\mathcal{F}_{\eta}$ is decreasing since $H^0(\mathcal{X},k\mathcal{L})$ is a $\mathbb{C}[t]$-module. We can extend $\mathcal{F}$ to a filtration by letting $$\mathcal{F}_{\eta}H^0(kL):=\mathcal{F}_{\lceil\eta\rceil}H^0(kL)$$ for non-integers $\eta,$ thus making $\mathcal{F}$ left-continuous. Since $$t^{-(\eta+\eta')}\overline{ss'}=(t^{-\eta}\bar{s})(t^{-\eta'}\bar{s'})\in H^0(\mathcal{X},k\mathcal{L})H^0(\mathcal{X},m\mathcal{L})\subseteq H^0(\mathcal{X},(k+m)\mathcal{L})$$ whenever $s\in \mathcal{F}_{\eta}H^0(kL)$ and $s'\in \mathcal{F}_{\eta'}H^0(kL),$ we see that $$(\mathcal{F}_{\eta} H^0(kL))(\mathcal{F}_{\eta'} H^0(mL))\subseteq \mathcal{F}_{\eta+\eta'}H^0((k+m)L),$$ i.e. $\mathcal{F}$ is multiplicative. Furthermore, by Lemma \ref{igenlemma} it follows that $\mathcal{F}$ is left-bounded and right-bounded. 

\begin{proposition} \label{hejho}
For $k\gg 0$ $$\mu(\mathcal{T},k)=\frac{d}{d\eta}(-\dim \mathcal{F}_{\eta}H^0(kL)).$$
\end{proposition}

\begin{proof}
Recall that we had the decomposition in weight spaces $$H^0(kL_0)=\oplus_{\eta}V_{\eta},$$ and that $$\mu(\mathcal{T},k):=\sum_{\eta=-\infty}^{\infty}\dim V_{\eta}\delta_{\eta}.$$ We have the following isomorphism: 
\begin{equation*}
(\pi_*k\mathcal{L})_{|\{0\}}\cong H^0(\mathcal{X},k\mathcal{L})/tH^0(\mathcal{X},k\mathcal{L}),
\end{equation*}
the right-to-left arrow being given by the restriction map, see e.g. \cite{Ross2}. Also, for $k\gg 0$, $(\pi_*k\mathcal{L})_{|\{0\}}=H^0(kL_0),$ therefore we get that for large $k$
\begin{equation} \label{isohej}
H^0(kL_0)\cong H^0(\mathcal{X},k\mathcal{L})/tH^0(\mathcal{X},k\mathcal{L}),
\end{equation}
We also had a decomposition of $H^0(\mathcal{X},k\mathcal{L})$ into the sum of its invariant weight spaces $W_{\eta}$. By Lemma \ref{igenlemma} it is clear that a section $S\in H^0(\mathcal{X},k\mathcal{L})$ lies in $W_{\eta}$ if and only if it can be written as $t^{-\eta}\bar{s}$ for some $s\in H^0(kL),$ in fact we have that $s=S_{|X}$. Thus we get that $$W_{\eta}\cong \mathcal{F}_{\eta}H^0(kL),$$ and by the isomorphism (\ref{isohej}) then $$V_{\eta}\cong \mathcal{F}_{\eta}H^0(kL)/\mathcal{F}_{\eta+1}H^0(kL).$$ Thus we get 
\begin{equation} \label{slutnej}
\dim \mathcal{F}_{\eta}H^0(kL)=\sum_{\eta'\geq \eta}\dim V_{\eta'},
\end{equation}
and the lemma follows by differentiating with respect to $\eta$ on both sides of the equation (\ref{slutnej}).
\end{proof}

\begin{proposition} \label{myprp}
The filtration associated to a test configuration $\mathcal{T}$ is always admissible. If we let $G_k[\mathcal{T}]$ denote the functions on $\Delta_k(L)$ associated to the filtration $\mathcal{F}(\mathcal{T})$ as previously definied, then we have that 
\begin{equation} \label{tatata}
\mu(\mathcal{T},k)=G_k[\mathcal{T}]_*\nu_k(L)
\end{equation}
and 
\begin{equation} \label{scaling}
\tilde{\mu}(\mathcal{T},k)=\frac{1}{k^n}((G_k[\mathcal{T}]/k)_*(\nu_k(L))).
\end{equation}
\end{proposition}

\begin{proof}
The equality of measures (\ref{tatata}) follows immediately from combining Lemma \ref{rara} and Proposition \ref{hejho}, and (\ref{scaling}) is just a rescaling of (\ref{tatata}). Since by Lemma \ref{phong} the weights of a test configuration is linearly bounded, by (\ref{tatata}) we get that the same holds for the functions $G_k[\mathcal{T}],$ i.e. the filtration $\mathcal{F}$ is linearly left- and right-bounded. It is hence admissible, since the other defining properties had already been checked.
\end{proof}

\begin{theorem} \label{mycor}
With the setting as in the proposition above, we have the following weak convergence of measures as $k$ tends to infinity $$\tilde{\mu}(\mathcal{T},k) \rightarrow G[\mathcal{T}]_*d\lambda_{|\Delta(L)}.$$
\end{theorem}

\begin{proof}
Follows from Theorem \ref{sebthm} together with Proposition \ref{myprp}.
\end{proof}

\begin{corollary}
In the asymtotic expansion $$\frac{w_k}{kd_k}=F_0-k^{-1}F_1+O(k^{-2})$$ we have that $$F_0=\frac{n!}{vol(L)}\int_{\Delta(L)}G(\mathcal{T})d\lambda.$$
\end{corollary}

\begin{proof}
Recall that in Section \ref{sec3} we defined $w_k$ by $$w_k:=\int_{\mathbb{R}} x d\mu(\mathcal{T},k),$$ i.e. in other words $$w_k=\sum \eta \dim V_{\eta},$$ $\oplus_{\eta}V_{\eta}$ being the weight space decomposition of $H^0(kL_0).$Thus Theorem \ref{mycor} implies that 
\begin{equation} \label{finally}
\lim_{k \to \infty}\frac{w_k}{k^{n+1}}=\lim_{k \to \infty}\int_{\mathbb{R}}x \tilde{\mu}(\mathcal{T},k)=\int_{\mathbb{R}}x(G[\mathcal{T}])_*(d\lambda_{|\Delta(L)})=\int_{\Delta(L)}G(\mathcal{T})d\lambda,
\end{equation}
using the weak convergence and the definition of the push forward of a measure. (\ref{finally}) together with the standard expansion $$d_k:=\dim H^0(kL)=k^n vol(L)/n!+o(k^n)$$ yields the corollary.
\end{proof}

Another consequence of Theorem \ref{mycor} is that it relates the Okounkov body $\Delta(L)$ with the central fibre $X_0,$ and therefore $X,$ in the sense of the following corollary.

\begin{corollary} \label{corowe}
Assume that we have embedded the test configuration $\mathcal{T}$ in some $\mathbb{P}^N\times \mathbb{C},$ let $h$ denote the corresponding Hamiltonian and $d\mu_{FS}$ the Fubini-Study volume measure on $|X_0|$ as in Section \ref{sec3}. Then we have that $$G[\mathcal{T}]_*d\lambda_{|\Delta(L)}=h_* d\mu_{FS}.$$ 
\end{corollary}

\begin{proof}
Follows immediately from combining Proposition \ref{reta} and Theorem \ref{mycor}.
\end{proof}

As in Section \ref{secny}, if restrict to the case of product test configurations where the $S^1$-action is effective, we can apply the Duistermaat-Heckman theorem to these measures, and get the following.

\begin{corollary} \label{retik}
Assume that there is a $\mathcal{C}^{\times}$-action on $X$ which lifts to $L,$ and that the corresponding $S^1$-action is effective. If we denote the associated product test configuration by $\mathcal{T},$ the concave transform $G[\mathcal{T}]$ is such that the pushforward measure $G[\mathcal{T}]_*d\lambda_{|\Delta(L)}$ is piecewise polynomial.
\end{corollary} 

\begin{proof}
Follows from combining Proposition \ref{prodtest} and Corollary \ref{corowe}.
\end{proof}

\section{Toric test configurations} \label{sec5}

We will cite some basic facts of toric geometry, all of which can be found in the article \cite{Donaldson2} by Donaldson.  
Let $L_P \to X_P$ be a toric line bundle with corresponding polytope $P\subseteq \mathbb{R}^n.$ Thus for every $k$ there is a basis for $H^0(kL_P)$ such that there is a one-one correspondence between the basis elements and the integer lattice points of $kP.$ We write this as $$\alpha\in kP\cap \mathbb{Z}^n \leftrightarrow z^{\alpha}\in H^0(kL_P).$$ In \cite{Donaldson2} Donaldson describes the relationship between toric test configurations and the geometry of polytopes. Let $g$ be a positive concave rational piecewise affine function defined on $P.$ One may define a polytope $Q$ in $\mathbb{R}^{n+1}$ with $P$ as its base and the graph of $g$ as its roof, i.e. $$Q:=\{(x,y): x\in P, y\in [0,g(x)]\}.$$ That $g$ is rational means precisely that the polytope $Q$ is rational, i.e. it is the convex hull of a finite set of rational points in $\mathbb{R}^n.$ In fact, by scaling we can without loss of generality assume that $Q$ is integral, i.e. the convex hull of a finite set of integer points. Then by standard toric geometry this polytope $Q$ corresponds to a toric line bundle $L_Q $ over a toric variety $X_Q$ of dimension $n+1.$ We may write the correspondence between integer lattice points of $kQ$ and basis elements for $H^0(kL_Q)$ as 
\begin{equation} \label{corrija}
(\alpha,\eta)\in kQ\cap \mathbb{Z}^{n+1} \leftrightarrow t^{-\eta}z^{\alpha}\in H^0(kL_Q).
\end{equation} 
There is a natural $\mathbb{C}^{\times}$-action $\rho$ given by multiplication on the $t$-variable. We also get a projection $\pi$ of $X_Q$ down to $\mathbb{P}^1,$ by letting $$\pi(x):=\frac{t^{-\eta+1}z^{\alpha}(x)}{t^{-\eta}z^{\alpha}(x)}$$ for any $\eta, \alpha$ such that this is well defined. Donaldson shows in \cite{Donaldson2} that if one excludes $\pi^{-1}(\infty),$ then the triple $L_Q, \rho$ and $\pi$ is in fact a test configuration, so $\pi$ is flat and the fiber over $1$ of $(X_Q, L_Q)$ is isomorphic to $(X_P,L_P).$

It was shown by Lazarsfeld-Musta\c{t}\u{a} in \cite{Lazarsfeld}, Example 6.1, that if one choses the coordinates, or actually the flag of subvarieties, so that it is invariant under the torus action, the Okounkov body of a toric line bundle is equal to its defining polytope, up to translation. Thus we may assume that $P=\Delta(L_P)$ and $$v(z^{\alpha})=\alpha.$$ The invariant meromorphic extension of the section $z^{\alpha}\in H^0(kL_P)$ is $z^{\alpha}\in H^0(kL_Q)$, where we have identified $X_P$ with the fiber over $1.$ By our calculations in Section \ref{sec4}, equation (\ref{useful}), the weight of $t^{-\eta}z^{\alpha}$ is $\eta.$ Thus we see that $$G_k(\alpha)=\sup \{\eta:  t^{-\eta}z^{k\alpha}\in H^0(kL_Q)\}=kg(\alpha),$$ by the correspondence (\ref{corrija}) and the fact that $g$ is the defining equation for the roof of $Q$. We get that $G_k/k$ is equal to the function $g$ restricted to $\Delta_k(L),$ and thus by the convergence of $G_k/k$ to $G[\mathcal{T}],$ that $$G[\mathcal{T}]=g.$$

We see that our concave transform $G[\mathcal{T}]$ is a proper generalization of the well-known correspondence between test configurations and concave functions in toric geometry.

It is thus clear that, as was shown for product test configurations in Proposition \ref{retik}, for toric test configurations it holds that the pushforward measure $$G[\mathcal{T}]_*d\lambda_{|L_P}=g_*d\lambda_{|P}$$ is the sum of a piecewise polynomial measure and a multiple of a dirac measure, simply because $P$ is a polytope and $g$ is piecewise affine (the dirac measure part coming the top of the roof).

\section{Deformation to the normal cone} \label{sec6}

One interesting class of test configurations is the ones which arise as a deformation to the normal cone with respect to some subscheme. This is described in detail by Ross-Thomas in \cite{Ross} and \cite{Ross2}, and we will only give a brief outline here.

Let $Z$ be any proper subscheme of $X.$ Consider the blow up of $X \times \mathbb{C}$ along $Z\times \{0\},$ and denote it by $\mathcal{X}.$ Hence we get a projection $\pi$ to $\mathbb{C}$ by composition $\mathcal{X} \to X \times \mathbb{C} \to \mathbb{C}.$ We let $P$ denote the exceptional divisor, and for any positive rational number $c$ we get a line bundle $$\mathcal{L}_c:=\pi^*L-cP.$$ By Kleimans criteria (see e.g. \cite{Lazarsfeld2}) it follows that $\mathcal{L}_c$ is relatively ample for small $c.$ The action on $(X \times \mathbb{C},L\times \mathbb{C})$ given by multiplication on the $\mathbb{C}$-coordinate lifts to an action $\rho$ on $(\mathcal{X}, \mathcal{L}_c),$ since both $Z\times \{0\}$ and $L\times \mathbb{C}$ are invariant under the action downstairs. Ross-Thomas in \cite{Ross} show that this data defines a test configuration.

From the proof of Theorem 4.2 in \cite{Ross} we get that 
\begin{equation} \label{billy}
H^0(\mathcal{X},k\mathcal{L}_c)=\bigoplus_{i=1}^{ck}t^{ck-i}H^0(X,kL\otimes \mathcal{J}_Z^{i})\oplus t^{ck}\mathbb{C}[t]H^0(kL),
\end{equation}
for $k$ sufficiently large and $ck\in \mathbb{N}.$ Here $\mathcal{J}_Z$ denotes the ideal sheaf of $Z,$ and the sections of $kL$ are being identified with their invariant extensions. From the expression (\ref{billy}) we can read off the associated filtration $\mathcal{F}$ of $H^0(kL).$ That $$t^{ck}H^0(kL)\subseteq H^0(\mathcal{X},k\mathcal{L}_c)$$ means that $$\mathcal{F}_{-ck}H^0(kL)=H^0(kL).$$ Furthermore, for $0 \leq i\leq ck$ and any $s\in H^0(kL)$ we get that $t^{ck-i}s\in H^0(\mathcal{X},k\mathcal{L}_c)$ iff $s\in H^0(kL\otimes \mathcal{J}_Z^{i}).$ This implies that for $-ck\leq \eta \leq 0,$ $$\mathcal{F}_{\eta}H^0(kL)=H^0(kL\otimes \mathcal{J}_Z^{ck+\eta}).$$ Also, when $\eta>0$ we get that $\mathcal{F}_{\eta}H^0(kL)=\{0\}.$ In summary, if we let $g_{c,k}$ be defined by $$g_{c,k}(\eta):=\lceil \max(\eta+ck,0)\rceil$$ for $\eta\in (-\infty,0]$ and let $g_{c,k}\equiv \infty$ on $(0,\infty),$ then by our calculations 
\begin{equation} \label{rarera}
\mathcal{F_{\eta}}H^0(kL)=H^0(kL\otimes \mathcal{J}_Z^{g_{c,k}(\eta)}).
\end{equation}
Thus this natural class of filtrations can be seen as coming from test configurations.

Let us assume that $Z$ is an ample divisor with a defining holomorphic section $s\in H^0(Z),$ i.e. $Z=\{s=0\}.$ Let $a$ be a number between zero and $c,$ then $L-aZ$ is still ample. Using multiplication with $s^{ka}$ we can embed $H^0(k(L-aZ))$ into $H^0(kL).$ With respect to this identification of $H^0(k(L-aZ))$ as a subspace of $H^0(kL)$ for all $k,$ we can identify the Okounkov body of $L-aZ$ with a subset of $\Delta(L).$ By vanishing theorems (see e.g. \cite{Lazarsfeld}), for large $k$
\begin{equation} \label{qwerty}
H^0(k(L-aZ))=H^0(kL\otimes \mathcal{J}_Z^{ka}),
\end{equation}
and therefore by (\ref{rarera})
$$H^0(k(L-aZ))=\mathcal{F}_{k(a-c)}H^0(kL).$$
It follows that the part of $\Delta(L)$ where $G[\mathcal{T}]$ is greater or equal to $a-c$ coincides with $\Delta(L-aZ).$\footnote{We thank Julius Ross for poining this out to us.}
 
Recall that by Theorem \ref{Lathm} $$\textrm{vol}_{\mathbb{R}^n}\Delta(L-aZ)=\frac{\textrm{vol}(L-aZ)}{n!}.$$ By this, a direct calculation yields that the pushforward measure $G[\mathcal{T}]_*d\lambda_{|\Delta(L)}$ can be written as $$\frac{\textrm{vol}(L-cZ)}{n!}\delta_0-\chi_{[-c,0]}\frac{d}{dx}\left(\frac{\textrm{vol}(L-(x+c)Z)}{n!}\right)dx,$$ where $\delta_0$ denotes the dirac measure at zero and $\chi_{[-c,0]}$ the indicator function of the interval $[-c,0].$ Since for any ample (or even nef)  ample line bundle the volume is given by integration of the top power of the first Chern class, $$\textrm{vol}(L)=\int_X c_1(L)^n,$$ it follows that the volume function is polynomial of degree $n$ in the ample cone. Thus the measure $G[\mathcal{T}]_*d\lambda_{|\Delta(L)}$ is a sum of a polynomial measure of degree less than $n$ and a dirac measure.

Let again $Z$ be an arbitrary subscheme of $X.$ Consider the blow up of $X$ along $Z,$ and let $E$ denote the exceptional divisor. If $E$ is irreducible we may introduce local holomorphic coordinates $(z_i)$ on the blow up, such that locally $E$ is given by the equation $z_1=0.$ Using these coordinates we get an associated Okounkov body $\Delta(L).$ For $s\in H^0(kL),$ the first coordinate of $v(s)$ is equal to the vanishing order of $s$ along $Z,$ i.e. the largest integer $r$ such that $s\in H^0(kL\otimes \mathcal{J}_Z^{r})$. Thus by (\ref{rarera}) we get that $$\Delta_{k,\eta}(L)=\{v(s)/k: s\in \mathcal{F}_{\eta}H^0(kL)\}=\Delta_k(L)\cap \{x_1\geq g_{c,k}(\eta)/k\}.$$ Furthermore 
\begin{eqnarray*}
G_k(\alpha)=\sup\{\eta: \alpha\in \Delta_{k,\eta}(L)\}=\\=\sup\{\eta: \alpha_1\geq g_{c,k}(\eta)/k\}=k\min(\alpha_1-c,0),
\end{eqnarray*}
and therefore $$G[\mathcal{T}](x)=\min(x_1-c,0).$$ 

\section{Product test configurations and geodesic rays} \label{sec7}

There is an interesting interplay between on the one hand test configurations and geodesic rays in the space of metrics on the other (see e.g. \cite{Sturm} and \cite{Sturm2}). The model case is when we have a product test configuration. 

Let $\mathcal{H}_L$ denote the space of positive hermitian metrics $\psi$ of a positive line bundle $L$ over $X.$ The tangent space of $\mathcal{H}_L$ at any point $\psi$ is naturally identified with the space of smooth real-valued functions on $X.$ The works of Mabuchi, Semmes and Donaldson (see \cite{Mabuchi}, \cite{Semmes} and \cite{Donaldson3}) have shown that there is a natural Riemannian metric on $\mathcal{H}_L,$ by letting the norm of a tangent vector $u$ at a point $\psi\in \mathcal{H}_L$ be defined by $$||u||_{\psi}^2:=\int_X |u|^2 dV_{\psi},$$ where $dV_{\psi}:=(dd^c\psi)^n.$ Let $\psi_t$ be a ray of metrics, $t\in (0,\infty)$. We may extend it to complex valued $t$ in $\mathbb{C}^{\times}$ if we let $\psi_t$ be independent on the argument of $t.$ We say that $\psi_t$ is a geodesic ray if 
\begin{equation} \label{geodesic}
(dd^c\psi_t)^{n+1}=0
\end{equation} 
on $X\times \mathbb{C}^{\times}.$ The equation (\ref{geodesic}) is the geodesic equation with respect to the Riemannian metric on $\mathcal{H}_L$ (see e.g. \cite{Sturm2}).

Let $\mathcal{T}$ be a product test configuration. That means that there is a $\mathbb{C}^{\times}$-action $\rho$ on the original pair $(X,L).$ Restriction of $\rho$ to the unit circle gives a $S^1$-action. Let $\varphi$ be an $S^1$-invariant positive metric on $L.$ We  get a $\mathbb{C}^{\times}$ ray $\tau \longmapsto \varphi_{\tau}\in \mathcal{H}_L$ of metrics by letting for any $\xi\in L$ 
\begin{equation} \label{mtric}
|\xi|_{\varphi_{\tau}}:=|\rho(\tau)^{-1}\xi|_{\varphi}.
\end{equation}
Similarly we get corresponding rays $k\varphi_{\tau}$ in $\mathcal{H}_{kL}.$
Since $\varphi$ was assumed to be $S^1$-invariant, $\varphi_{\tau}$ only depends on the absolute value $|\tau|.$ Also because the action $\rho$ is holomorphic, it follows that $$(dd^c\varphi_{\tau})^{n+1}=0,$$ therefore $\varphi_{\tau}$ is a geodesic ray. 

In \cite{Berndtsson} Berndtsson introduces sequences of spectral measures on $\mathbb{R}$ arising naturally from a geodesic segment of metrics, and shows that they converge weakly to a certain pushforward of a volume form on $X$. Inspired by his result, we consider the analogue in our setting. 

Let $\dot{\varphi}$ denote the derivative of $\varphi_{\tau}$ at $1,$ so $\dot{\varphi}$ is a smooth real-valued function on $X.$ We consider the positive measure on $\mathbb{R}$ we get by pushing forward the volume form $dV_{\varphi}:=(dd^c\varphi)^n$ on $X$ with this function divided by two, $$\mu_{\varphi}:=(\dot{\varphi}/2)_*dV_{\varphi}.$$ The measure $\mu_{\varphi}$ does not does not depend on the choice of $S^1$-invariant metric $\varphi.$ In fact, we have the following result.

\begin{theorem} \label{bigshot}
Let $G[\mathcal{T}]$ denote the concave transform of the product test configuration. We have an equality of measures $$\mu_{\varphi}=G[\mathcal{T}]_*d\lambda_{|\Delta(L)}.$$
\end{theorem}

\begin{proof}
We will use one of the main ideas in the proof of the main result of Berndtsson in \cite{Berndtsson}, Theorem 3.3. However, in our setting where the geodesic comes from a $\mathbb{C}^{\times}$-action things are much simpler since we do not need the powerful estimates used in \cite{Berndtsson}. 

Let $dV$ be some fixed smooth volume form on $X.$ We will introduce two families of scalar products on $H^0(kL),$ parametrized by $\tau,$ $||.||_{\tau,1}$ and $||.||_{\tau,2}.$ First we let for any $s\in H^0(kL)$ $$||s||^2_{\tau,1}:=\int_X |s|^2_{k\varphi_{\tau}}dV,$$ 
while we let $$||s||^2_{\tau,2}:=\int_X |\rho(\tau)^{-1}s|^2_{k\varphi}dV= ||\rho(\tau)^{-1}s||_{1,1}^2.$$ 

Direct calculations yield that 
\begin{equation} \label{shortline}
\frac{d}{d\tau}||s||^2_{\tau,1}=\frac{d}{d\tau}\int_X |s|^2_{k\varphi_{\tau}}dV=\int_X (-k\dot{\varphi}_{\tau})|s|^2_{k\varphi_{\tau}}dV=(T_{-k\dot{\varphi}_{\tau}}s,s)_{\tau,1},
\end{equation}
where $T_{-k\dot{\varphi}_{\tau}}$ denotes the Toeplitz operator with symbol $-k\dot{\varphi}_{\tau}.$ 

Differentiating $||.||_{\tau,2}$ with respect to $\tau$ we get that 
\begin{equation} \label{shortline2}
\frac{d}{d\tau}||s||^2_{\tau,2}=\frac{d}{d\tau}(\rho(\tau)^{-1}s,\rho(\tau)^{-1}s)_{1,1}=((\frac{d}{d\tau}\rho(\tau)^{-2})s,s)_{1,1}.
\end{equation}

On the other hand
\begin{eqnarray} \label{longline}
||s||^2_{\tau,1}=\int_X |s(x)|^2_{k\varphi_{\tau}}dV(x)=\int_X |\rho(\tau)^{-1}(s(x))|^2_{k\varphi}dV(x)= \nonumber \\=\int_X |(\rho(\tau)^{-1}s)(x)|^2_{k\varphi}dV(\rho(\tau)x)=\int_X |\rho(\tau)^{-1}s|^2_{k\varphi}dV_{\tau},
\end{eqnarray} where $dV_{\tau}(x):=dV(\rho(\tau)x)$ thus denotes the resulting volume form after the $\tau$-action. Since $dV_{\tau}(x)$ depends smoothly on $\tau,$ using (\ref{longline}) we get that 
\begin{eqnarray} \label{longline2}
\left |\frac{d}{d\tau}_{|\tau=1}||s||^2_{\tau,1}-\frac{d}{d\tau}_{|\tau=1}||s||^2_{\tau,2}\right |=\left |\frac{d}{d\tau}_{|\tau=1}\int_X |\rho(\tau)^{-1}s|^2_{k\varphi}(dV_{\tau}-dV)\right |\leq \nonumber \\ \leq \int_X |\frac{d}{d\tau}_{|\tau=1}dV_{\tau}| \int_X |s|^2_{k\varphi}dV=C||s||^2_{1,1},
\end{eqnarray} 
where thus $C$ is a uniform constant independent of $s$ and $k.$ 
Therefore letting $\tau=1$ in equations (\ref{shortline}) and (\ref{shortline2}), and using (\ref{longline2}) we get that 
\begin{equation} \label{import}
\frac{d}{d\tau}\rho(\tau)_{|\tau=1}=T_{k\dot{\varphi}/2}+E_k,
\end{equation}
where the error term $E_k$ is uniformly bounded, $||E_k||<C'.$

Let $A$ be a self-adjoint operator on a $N$-dimensional Hilbert space, and let $\lambda_i$ denote the eigenvalues of $A,$ which therefore are real, counted with multiplicity. The spectral measure of $A,$ denoted by $\nu(A),$ is defined as $$\nu(A):=\sum_i \delta_{\lambda_i}.$$ 

We consider the normalized spectral measure of $T_{k\dot{\varphi}/2}$, $$\nu_k:=\frac{1}{k^n}\nu(T_{k\dot{\varphi}/2}/k).$$ By Theorem 3.2 in \cite{Berndtsson}, which is a variant of a theorem of Boutet de Monvel-Guillemin (see \cite{Boutet}), we get that the measures $\nu_k$ converge weakly as $k$ tends to infinity to the measure $\mu_{\varphi}.$ 

Let $H^0(kL)=\sum_{\eta}V_{\eta}$ be the decomposition in weight spaces, and let $P_{\eta}$ denote the projection to $V_{\eta}.$ Then $$\rho(\tau)=\sum_{\eta}\tau^{\eta}P_{\eta},$$ and thus
\begin{equation} \label{opeq2}
\frac{d}{d\tau}\rho(\tau)_{|\tau=1}=\sum \eta P_{\eta}.
\end{equation}
From (\ref{opeq2}) we see that the normalized spectral measures of $\frac{d}{d\tau}\rho(\tau)_{|\tau=1},$ which we denote by $\mu_k,$ coincides with the previously defined weight measure $$\tilde{\mu}(\mathcal{T},k)=\frac{1}{k^n}\sum_{\eta=-\infty}^{\infty}\dim V_{\eta}\delta_{k^{-1}\eta}.$$ According to Theorem \ref{mycor} the sequence $\tilde{\mu}(\mathcal{T},k)$, and therefore $\mu_k,$ converges weakly to the measure $G[\mathcal{T}]_*d\lambda_{|\Delta(L)}$. 

Lastly, by the the min-max principle, when perturbing an operator $A$ by an operator $E$ with small norm $||E||<\varepsilon,$ then each eigenvalue is perturbed at most by $\varepsilon.$ Thus from (\ref{import}) it follows that $\nu_k-\mu_k$ converges weakly to zero, and the theorem follows.
\end{proof}

We will relate this result to our previous discussion on Duistermaat-Heckman measures in Section \ref{secny} and \ref{sec4}, by showing that the map $\dot{\varphi}/2$ is a Hamiltonian for the $S^1$-action when the symplectic form is given by $dd^c\varphi$. This is of course well-known (see e.g. \cite{Donaldson3}), but we include it here for the benefit of the reader.

Let $V$ be the holomorphic vector field on $X$ generating the action $\rho.$ Hence, the imaginary part $\textrm{Im}V$ of $V$ generates the $S^1$-action. By definition, $\dot{\varphi}/2$ is a Hamiltonian if it holds that 
\begin{equation} \label{trewu}
\textrm{Im}V\rfloor dd^c\varphi=d\dot{\varphi}/2,
\end{equation} 
where $\rfloor$ denotes the contraction operator.

If we can show that $$-iV\rfloor dd^c\varphi=\bar{\partial}\dot{\varphi}/2,$$ equation (\ref{trewu}) will follow by taking the real part on both sides. We calculate locally with respect to some trivialization and without loss of generality we may assume that $$V=\frac{\partial}{\partial z_1}.$$ Recall that by definition $$dd^c\varphi=\frac{i}{2}\sum \frac{\partial^2\varphi}{\partial z_i \partial \bar{z}_j}dz_i\wedge d\bar{z}_j.$$ Hence we get that $$-iV\rfloor dd^c\varphi=\frac{1}{2}\sum \frac{\partial^2\varphi}{\partial z_1 \partial \bar{z}_j}d\bar{z}_j=\frac{1}{2}\bar{\partial}\frac{\partial \varphi}{\partial z_1}.$$ Since $V=\partial/\partial z_1$ generates the action, it follows that locally $\partial/\partial z_1\varphi=\dot{\varphi},$ and we are done.

\end{document}